\documentclass[a4paper]{cdc}
\usepackage{amstex}

\usepackage{epsfig}
\def\figum{deg_1.eps}
\def\figdois{deg_2.eps}
\def\figtres{sin_1.eps}
\def\figquatro{sin_2.eps}
\def\figcinco{sia_1.eps}
\def\figseis{sia_2.eps}
\def\figsete{prb_1.eps}
\def\figoito{prb_2.eps}

\def\rep#1{(\ref{#1})}
\def\FAPESP{\textsc{fapesp}}
\def\CNPq{\textsc{cnp}q}

\begin{document}

\title{A Tuner that \\ Accelerates Parameters}
\author{{Felipe M Pait\thanks{Research supported  by {\CNPq} -- Brazilian Research Council, grant 523949/95-2.}
 and Paulo A Atkinson\thanks{Research supported  by {\FAPESP} -- S\~ao Paulo Research Agency, scholarship 97/00582-5.}} \\ 
Universidade de S\~ao Paulo \\
Laborat\'orio de Automa\c c\~ao e Controle  -- \textsc{pee}\\ 
S\~ao Paulo SP 05508--900 Brasil \\
{\tt pait,atk\char'100lac.usp.br} }
\maketitle
\thispagestyle{empty}\pagestyle{empty}

\begin{abstract} We propose a tuner, suitable for adaptive control and (in its discrete-time version) adaptive filtering applications, that sets the second derivative of the parameter estimates rather than the first derivative as is done in the overwhelming majority of the literature. Comparative stability and performance analyses are presented.
\vspace{1ex}
 
 \end{abstract}

\section{Parameter Adjustment}In adaptive control and recursive parameter estimation one often needs to adjust recursively an estimate $\hat{p}$ of a vector $p$, which comprises $n$ constant but unknown parameters, using measurements of a quantity
\begin{equation}\label{eq:y}
y = x^\top p  + w.
\end{equation}
Here $x : [0,\bar{t}) \rightarrow  \Bbb{R}^n$ is a vector of known data, often called the regressor, and $w$ is a measurement error signal. 
The goal of tuning is to keep both the estimation error $x^\top \hat{p} - y$ and the parameter error $\hat{p} - p$ as small as possible.

There are several popular methods for dealing with the problem above, for instance least-squares. Maybe the most straightforward involve minimizing the prediction error via gradient-type algorithms of the form:
\begin{equation}\label{eq:phat}
\Dot{\Hat{p}} = -Mx (x^\top \hat{p} - y),
\end{equation}
where $M$ is a constant, symmetric, positive-definite gain matrix. Let us define $q = \hat{p} - p$ and analyze differential equations \rep{eq:y} and \rep{eq:phat}, which under the assumption that $w(t)$ is identically zero read:
\begin{equation}\label{eq:qdot}
\dot{q} = - Mx x^\top q.
\end{equation}
The nonnegative function 
$V = \frac{1}{2} q^\top M^{-1} q$
has time derivative 
$\dot{V} = q^\top M^{-1} \dot{q} = -q^\top x x^\top q,$ hence 
\[\int_0^t \dot{V} = V(t) - V(0) = -\int_0^t  (x^\top q)^2.\]

Inspection of the equation above reveals that $V$ is limited in time, thus $q \in {\mathcal L}^\infty$, and also that the error $x^\top q  \in {\mathcal L}^2$ (norms are taken on the interval $[0,\bar{t})$ where all signals are defined). These are the main properties an algorithm needs in order to be considered a suitable candidate for the role of a tuner in an adaptive control system. Often  $\dot{q}  \in {\mathcal L}^2$ or something similar is also a  desirable property. To obtain the latter,  normalized algorithms can be used; however, the relative merits of normalized versus unnormalized tuners are still somewhat controversial.  Another alternative is to use a time-varying $M$, as is done  in least-squares tuning.

In \S\ref{sec:acceleration} we present a tuner that sets the second derivative of $\hat{p}$, and in \S\ref{sec:covariance} the effects of a white noise $w$ on the performance of the two algorithms are compared.  Then we show some simulations and make concluding remarks.

\section{The Accelerating Algorithm} \label{sec:acceleration}
Classical tuners are such that the \textsl{velocity\/} of adaptation (the first derivative of the parameters) is set proportional to the regressor and to the prediction error $x^\top \hat{p} - y = x^\top q$. We propose  to set the \textsl{acceleration\/} of the parameters:
\begin{equation}\label{eq:qddot}
\boxed
{\ddot{q} = - xx^\top q - 2 (I + xx^\top)\dot{q}.}
\end{equation}
Notice that the the formula above is implementable (using $2n$ integrators) if measurement error is absent, because the unknown $q$ appears only in scalar product with $x$. Choose another function of Lyapunovian inspiration: 
\begin{align*}
V &= q^\top q + q^\top \dot{q} + \dot{q}^\top \dot{q} \\
    &= \begin{bmatrix} q^\top & \dot{q}^\top \end{bmatrix}
\begin{bmatrix} I & I/2 \\ I/2 & I \cr \end{bmatrix}
\begin{bmatrix} q \\ \dot{q} \\ \end{bmatrix} \geq 0.
\end{align*}
Taking derivatives along the trajectories of  \rep{eq:qddot} gives
\begin{align*}
\dot{V} &= 
2q^\top \dot{q} + \dot{q}^\top \dot{q} 
        - (q^\top + 2\dot{q}^\top) (2 \dot{q} + xx^\top q + 2xx^\top \dot{q}) \\
             &=  - 3 \dot{q}^\top \dot{q} -  (q + 2\dot{q})^\top xx^\top (q + 2\dot{q})  \leq 0.
\end{align*}
Integrating $\dot{V}$ we obtain 
\[
V(t) = V(0) -\int_0^t  \dot{q}^\top \dot{q} - \int_0^t \Bigl(x^\top (q+ 2\dot{q})\Bigr)^2,
\]
which leads immediately to the desired properties:
\[
q,\dot{q} \in \mathcal{L}^\infty; \quad \dot{q} \in \mathcal{L}^2; \quad x^\top(q+2\dot{q}) \in \mathcal{L}^2.
\]

The slow variation property $\dot{q} \in \mathcal{L}^2$ follows without the need for normalization, and now we obtain $x^\top (q+2\dot{q}) \in \mathcal{L}^2$ instead of $x^\top q \in \mathcal{L}^2$ as before. We might regard $x^\top (q+2\dot{q}) $ as a modified error, which can be used in the stability analysis of a detectable or ``tunable'' adaptive system via an output-injection argument; see \cite{towards2}. A generalization of \rep{eq:qddot} is
\begin{equation}\label{eq:qddot2}
\ddot{q} = - M_1 \Bigl ( x(x^\top \hat{p} -y) + 2 (M_2 + xx^\top M_1 M_3) \dot{q} \Bigr ),
\end{equation}
with $M_1, M_2$ and $M_3$ constant, symmetric, positive-definite $n\times n$ matrices such that $M_2^{-1} < 4 M_1 M_3 M_1$ and $M_2M_1M_3 + M_3 M_1 M_2 > M_1^{-1} /2$. The properties of tuner  \rep{eq:qddot2}, which can be obtained using the positive-definite function
\[
\begin{bmatrix}q^\top & \dot{q}^\top \end{bmatrix}
\begin{bmatrix} M_2 & M_1^{-1}/2 \\ M_1^{-1}/2 & M_3  \end{bmatrix}
\begin{bmatrix} q \cr \dot{q}\end{bmatrix}
\]
in the same manner as before, are
\[
q,\dot{q} \in \mathcal{L}^\infty; \quad \dot{q} \in \mathcal{L}^2; \quad x^\top(q+2M_1 M_3 \dot{q}) \in \mathcal{L}^2.
\]

\section{Covariance Analysis} \label{sec:covariance}

We now consider the effects on the expected value and covariance of $q(t)$ of the presence of a measurement error. The assumptions are that $w(t)$ is a white noise with zero average and covariance $\sigma_w^2 \delta(t-\tau)$ and that $x(t)$ are given, deterministic data. For comparison purposes, first consider what happens when the conventional tuner \rep{eq:phat}  is applied to \rep{eq:y} in the presence of measurement error $w$:
\begin{equation}\tag{\ref{eq:qdot}$'$} \label{eq:qdot'}
\dot{q} = - Mx x^\top q + Mx w.
\end{equation}
The solution to the equation above can be written in terms of $(-Mxx^\top)$'s state transition matrix $\phi(t,\tau)$  as follows
\[
q(t) = \phi(t,0) q(0) + \int_0^t d\tau \, \phi(t,\tau) M x(\tau) w(\tau),
\]
hence $E_w\{q(t)\} = \phi(t,0) q(0)$ because $E\{w(\tau)\} = 0$ by assumption. Here the notation $E_w\{\cdot\}$, denoting the expectation with respect to the random variable $w$, is used to emphasize that the stochastic properties of $x$ are not under consideration. The conclusion is that $q(t)$ will converge to zero in average as fast as $\phi(t,0)$ does. The well-known persistency of excitation conditions on $x(\cdot)$ are sufficient for the latter to happen.

To study the second moment of the parameter error, write
\begin{equation*}
\begin{split}
q q^\top(t) & = \phi(t,0) qq^\top(0) \phi^\top(t,0) \\
	& \quad+ \phi(t,0) q(0) \int_0^t d\tau \, w(\tau) x^\top(\tau) M \phi^\top(t,\tau) \\
	& \quad+ \left(\int_0^t d\tau \, \phi(t,\tau) M x(\tau) w(\tau)\right) q^\top(0) \phi^\top(t,0) \\
	& \quad + \left(\int_0^t d\tau \, \phi(t,\tau) M x(\tau) w(\tau)\right) \cdot \\
	& \qquad \qquad \left( \int_0^t d\tau \, w(\tau) x^\top(\tau) M \phi^\top(t,\tau)\right). \\
\end{split}\end{equation*}
The covariance of $q(t)$ can be written as the sum of four terms.
The first is deterministic. 
%
The second term 
\begin{multline*}
E_w\Bigl\{\phi(t,0) q(0) \int_0^t d\tau \, w(\tau) x^\top(\tau) M \phi^\top(t,\tau)\Bigr\} \\
 = \phi(t,0) q(0) \int_0^t d\tau \, E\{w(\tau)\} x^\top(\tau) M \phi^\top(t,\tau) = 0
\end{multline*}
 because $w(t)$ has zero mean, and the third term is likewise zero.
The fourth term
\begin{equation*} \begin{split}
& E_w \left \{ \int_0^t d\tau_1 \, \phi(t,\tau_1) M x(\tau_1) w(\tau_1) \right. \\
& \left. \qquad \qquad \qquad \int_0^t d\tau_2 \, w(\tau_2) x^\top(\tau_2) M \phi^\top(t,\tau_2) \right \} \\
& \qquad = \int_0^t d\tau_1\int_0^t d\tau_2 \, \phi(t,\tau_1) M x(\tau_1)  \\
& \qquad \qquad \qquad \qquad E\{w(\tau_1) w(\tau_2)\} x^\top(\tau_2) M \phi^\top(t,\tau_2)  \\
& \qquad = \sigma_w^2 \int_0^t d\tau_1 \, \phi(t,\tau_1) M x(\tau_1)  x^\top(\tau_1) M \phi^\top(t,\tau_1), \\
\end{split} \end{equation*}
where  Fubini's Theorem and the fact $E\{w(\tau_1)w(\tau_2)\} = \sigma_w^2 \delta(\tau_1-\tau_2)$ were used. Performing the integration and adding the first and fourth terms  results in
\begin{multline*}
E_w\{qq^\top(t)\} = \phi(t,0) q q^\top(0) \phi^\top(t,0)  \\
+  \frac{\sigma_w^2}{2}\Bigl(M - \phi(t,0) M \phi^\top(t,0)\Bigr).
\end{multline*}
This equation can be given the following interpretation: for small $t$, when $\phi(t,0)$ is close to the identity, the covariance of $q(t)$ remains close to $qq^\top(0)$, the outer product of the error in the initial guess of the parameters with itself. As $\phi(t,0) \rightarrow 0$, which will happen if $x(\cdot)$ is persistently exciting, $E_w\{ qq^\top(t)\}$ tends to $\sigma_w^2 M /2$. This points to a compromise between higher convergence speeds and lower steady-state parameter error, which require respectively larger and smaller values of the gain $M$. Algorithms that try for the best of both worlds --- parameter convergence in the mean-square sense ---  often utilize time-varying, decreasing gains; an example is the least-squares algorithm.

We shall now attempt a similar analysis for the acceleration tuner \rep{eq:qddot2} applied to \rep{eq:y}, which results in the differential equation
\begin{equation*}
\begin{bmatrix} \dot{q} \\ \ddot{q} \end{bmatrix} = 
\begin{bmatrix} 0 & I \\ \scriptstyle -M_1xx^\top &\scriptstyle -2M_1(M_2+ xx^\top M_1M_3) \end{bmatrix} \begin{bmatrix} q \\ \dot{q} \end{bmatrix} + \begin{bmatrix} 0 \\ \scriptstyle M_1 x w \end{bmatrix}.
\end{equation*}
Let
\begin{equation*}
\begin{bmatrix} \dot{\phi}_{11}\!\! & \!\!\dot{\phi}_{12} \\ \dot{\phi}_{21}\!\! & \!\!\dot{\phi}_{22} \end{bmatrix} = 
\begin{bmatrix} 0 \!\!&\!\! I \\ \scriptstyle -M_1xx^\top \!\!&\!\! \scriptstyle -2M_1(M_2+ xx^\top M_1M_3) \end{bmatrix}  
\begin{bmatrix} {\phi}_{11} \!\!&\!\! {\phi}_{12} \\ {\phi}_{21} \!\!&\!\!{\phi}_{22} \end{bmatrix},
\end{equation*}
where  ${\phi}_{11}(t,t) = {\phi}_{22}(t,t) = I$, ${\phi}_{12}(t,t) = {\phi}_{21}(t,t) =0$,
each $\phi_{ii}$ is a function of $(t,\tau)$ unless otherwise noted, and the dot signifies derivative  with respect to the first argument. If  $\dot{q}(0) = 0$,
\begin{equation*}
q(t) = \phi_{11}(t,0) q(0) + \int_0^t d\tau \, \phi_{12}(t,\tau) M_1 x(\tau) w(\tau).
\end{equation*}

Following the same reasoning used for the velocity tuner, one concludes that  $E_w\{q(t)\} = \phi_{11}(t,0) q(0)$ and that
\begin{multline} \label{eq:cov-acelera}
E_w\{qq^\top(t)\} = \phi_{11}(t,0) q q^\top(0) \phi_{11}^\top(t,0) + \\
\sigma_w^2 \int_0^t d\tau \, \phi_{12}(t,\tau) M_1 x x^\top(\tau) M_1 \phi_{12}^\top(t,\tau),
\end{multline}
however the properties of the acceleration and velocity tuners are not yet directly comparable because  the right-hand side of \rep{eq:cov-acelera} does not lend itself to immediate integration. To obtain comparable results, we employ the ungainly but easily verifiable formula \rep{eq:integration-formula},
\begin{figure*}[hbt]
\begin{multline} \label{eq:integration-formula}
{\partial \over \partial \tau}{
\begin{bmatrix}{\phi}_{11} & {\phi}_{12} \\ {\phi}_{21} &{\phi}_{22} \end{bmatrix} \bar{M}_1
\begin{bmatrix}{\phi}^\top_{11} & {\phi}^\top_{21} \\ {\phi}^\top_{12} &{\phi}^\top_{22} \end{bmatrix}}  
= \left.\begin{bmatrix}{\phi}_{11} & {\phi}_{12} \\ {\phi}_{21} &{\phi}_{22} \end{bmatrix}  \right(
\bar{M}_2  \\
+ \left. \begin{bmatrix}0 & 0 \\ 0 & 
-M_1xx^\top M_1 + 2(1+\mu_2)  (M_1M_2M_1M_3M_1xx^\top M_1+M_1xx^\top M_1M_3M_1M_2M_1) 
 \end{bmatrix}
 \right)  \begin{bmatrix}{\phi}^\top_{11} & {\phi}^\top_{21} \\ {\phi}^\top_{12} &{\phi}^\top_{22} \end{bmatrix}, \quad  \\
\begin{split}
\text{where}\quad  \bar{M}_1 &= \begin{bmatrix}M_1 M_3 M_1 + \mu_1 M_1 & -M_1/2 \\ -M_1/2 &(1+\mu_2)M_1M_2M_1 \end{bmatrix}, \\
  \bar{M}_2 &= \begin{bmatrix}M_1  & \mu_1 M_1 xx^\top M_1 - (2+\mu_2) M_1 M_2 M_1  \\  \mu_1 M_1 xx^\top M_1-(2+\mu_2) M_1 M_2 M_1  &4(1+\mu_2)M_1M_2M_1M_2M_1 \end{bmatrix} \quad 
\end{split}
\end{multline}
\rule{\linewidth}{0.2pt}
\end{figure*}
valid for arbitrary scalars $\mu_1$ and $\mu_2$, and make the

\paragraph{Simplifying Assumption: }{} For $i =1,2$, and 3, $M_i = m_i I$, where $m_i$ are scalars and $I$ is the $n \times n$ identity matrix.

Premultiplying \rep{eq:integration-formula} by $[ I \; 0 ]$, 
postmultiplying  by $[ I \; 0 ]^\top$,   integrating from 0 to $t$, and using the simplifying assumption gives formula \rep{eq:integration-formula2}.
\begin{figure*}[t]
\begin{multline} \label{eq:integration-formula2}
\int_0^t d\tau \, (4(1+\mu_2)m_1^2 m_2m_3 -1) \phi_{12}(t,\tau) M_1 x  x^\top(\tau) M_1 \phi_{12}^\top(t,\tau) \\
=   \begin{bmatrix}{\phi}_{11}(t,\tau) & {\phi}_{12}(t,\tau)  \end{bmatrix} \bar{M}_1
\begin{bmatrix}{\phi}^\top_{11}(t,\tau) \\ {\phi}^\top_{12}(t,\tau)  \end{bmatrix}\Biggr |_{\tau = 0}^{\tau = t} 
 - \int_0^t d\tau \begin{bmatrix}{\phi}_{11} & {\phi}_{12}  \end{bmatrix} \bar{M}_{2}
\begin{bmatrix}{\phi}^\top_{11} \\ {\phi}^\top_{12}  \end{bmatrix}
\end{multline}
\rule{\linewidth}{0.2pt}
\end{figure*}
Taking $\mu_1 = \mu_2 = 0$ in \rep{eq:integration-formula2}, $\bar{M}_2$ results positive-semidefinite, therefore
\begin{multline} \label{eq:desigualdade}
\int_0^t d\tau \, \phi_{12}(t,\tau) M_1 x  x^\top(\tau) M_1 \phi_{12}^\top(t,\tau) \\
\leq   \frac{1}{4m_1^2 m_2m_3 -1} \left(m_1^2m_3 I 
\phantom{\begin{bmatrix}{\phi}^\top_{11}(t,0) \\ {\phi}^\top_{12}(t,0)  \end{bmatrix}}\right. \\
\left. - \begin{bmatrix}{\phi}_{11}(t,0) & {\phi}_{12}(t,0)  \end{bmatrix} \bar{M}_1
\begin{bmatrix}{\phi}^\top_{11}(t,0) \\ {\phi}^\top_{12}(t,0)  \end{bmatrix}\right).
\end{multline}

The combination of \rep{eq:cov-acelera} and \rep{eq:desigualdade} shows that $m_1$ can be increased without affecting  $q$'s steady-state covariance. On the other hand, to decrease the covariance we need to increase $m_2$, which roughly speaking means increasing damping in \rep{eq:qddot}. Since $m_1$ and $m_2$ can be increased  without affecting the stability properties shown in \S\ref{sec:acceleration}, a better transient $\times$ steady-state performance compromise might be achievable with the acceleration tuner than with the velocity tuner, at least in the case when $M_1$, $M_2$, and $M_3$ are ``scalars.''  Notice that $4m_1^2 m_2m_3 >1$ by construction.

\paragraph{Approximate  Analysis: }{} The derivation of inequality \rep{eq:desigualdade} does not involve any approximations, and therefore provides an upper bound on $E_w\{qq^\top(t)\}$, valid independently of $x(\cdot)$. A less conservative estimate of the integral in \rep{eq:cov-acelera} can be obtained by replacing $xx^\top$ by its average value $R$ in the definition of $\bar{M}_2$ in \rep{eq:integration-formula2}. This approximation seems reasonable because $\bar{M}_2$ appears inside an integral, but calls for more extensive simulation studies. 

To obtain a useful inequality, we require $\bar{M}_2 \geq 0$; namely, using the Schur complement
\begin{multline*}
4(1+\mu_2) M_1M_2M_1M_2M_1 \geq \\
 M_1(\mu_1 xx^\top-(2+\mu_2) M_2) M_1(\mu_1 xx^\top-(2+\mu_2) M_2) M_1,
\end{multline*}
or, using the simplifying assumption and substituting $xx^\top$ by its approximation $R$
\[
4(1+\mu_2)I \geq \Bigl( (2 +\mu_2)I - \frac{\mu_1 }{m_2}R  \Bigr)^2.
\]
Suppose further that $R = rI$. Looking for the least conservative estimate, we pick $\mu_1 = \frac{m_2}{r}(2+\mu_2 -2\sqrt{1+\mu_2})$, the least value of $\mu_1$ that keeps $\bar{M}_2 \geq 0$. Thus
\begin{multline} \label{eq:desigualdade-approx}
\int_0^t d\tau \, \phi_{12}(t,\tau) M_1 x  x^\top(\tau) M_1 \phi_{12}^\top(t,\tau) \leq\\
   \frac{m_1^2m_3 r +m_1 m_2 (\mu_2+2 -2\sqrt{1+\mu_2})}{4m_1^2 m_2m_3r(1+\mu_2) -r}  I - \Xi(t),
\end{multline}
with $\Xi(t) = 
\frac{\left[\begin{smallmatrix}{\phi}_{11}(t,0) & {\phi}_{12}(t,0)  
\end{smallmatrix}\right] \bar{M}_1
\left[\begin{smallmatrix}{\phi}^\top_{11}(t,0) \\ {\phi}^\top_{12}(t,0)  \end{smallmatrix}\right]}{4m_1^2 m_2m_3r(1+\mu_2) -r}.$

Taking $\mu_2 = 0$ we repeat the previous, exact result. For large positive values of $\mu_2$ the first term of the right-hand side of \rep{eq:desigualdade-approx} tends to $I / 4 m_1m_3r$, which indicates that the steady-state covariance of the parameter error decreases when the signal $x$ increases in magnitude, and that it can be made smaller via appropriate choices of the gains $m_1$ and $m_3$. The situation for the accelerating tuner is hence much more favorable than for the conventional one.

\section{Simulations}
\begin{figure*}[t] 
 \centerline{\epsfysize=2.5in \epsffile{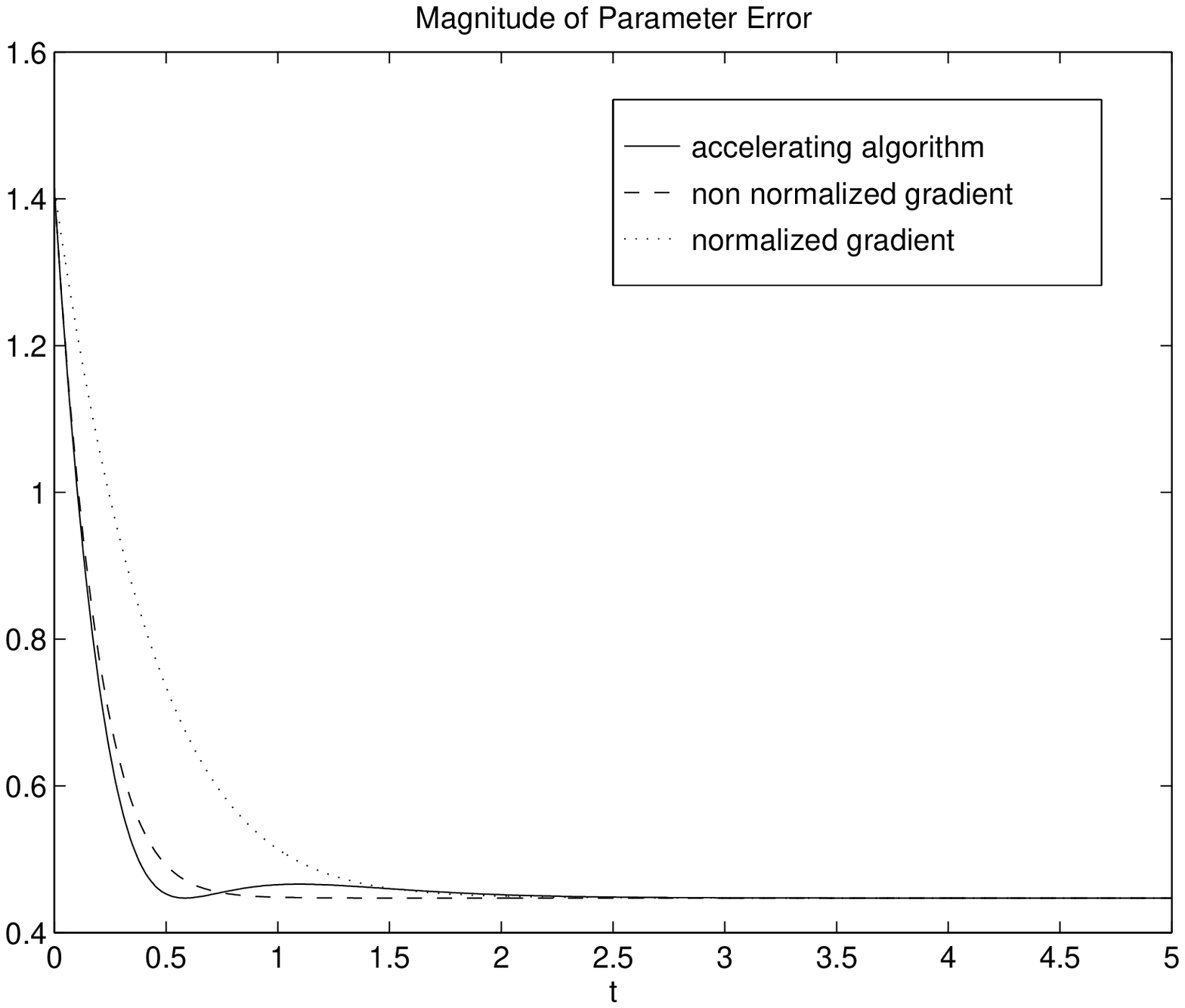} \epsfysize=2.5in \epsffile{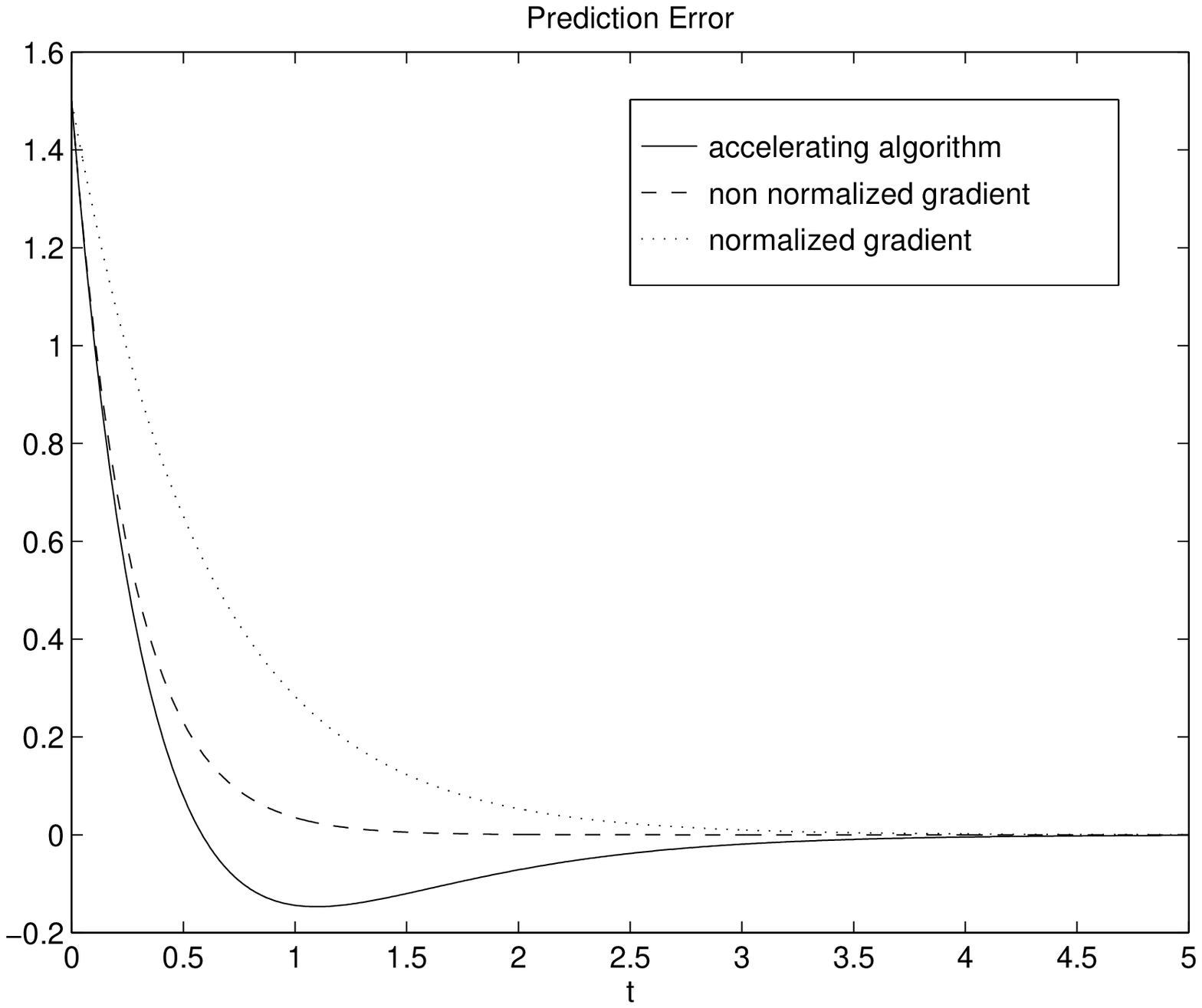}} 
 \caption{Response to step data}
 \label{fig:step}
\end{figure*}

\begin{figure*}[p!] 
 \epsfysize=2in
 \centerline{\epsfysize=2.5in \epsffile{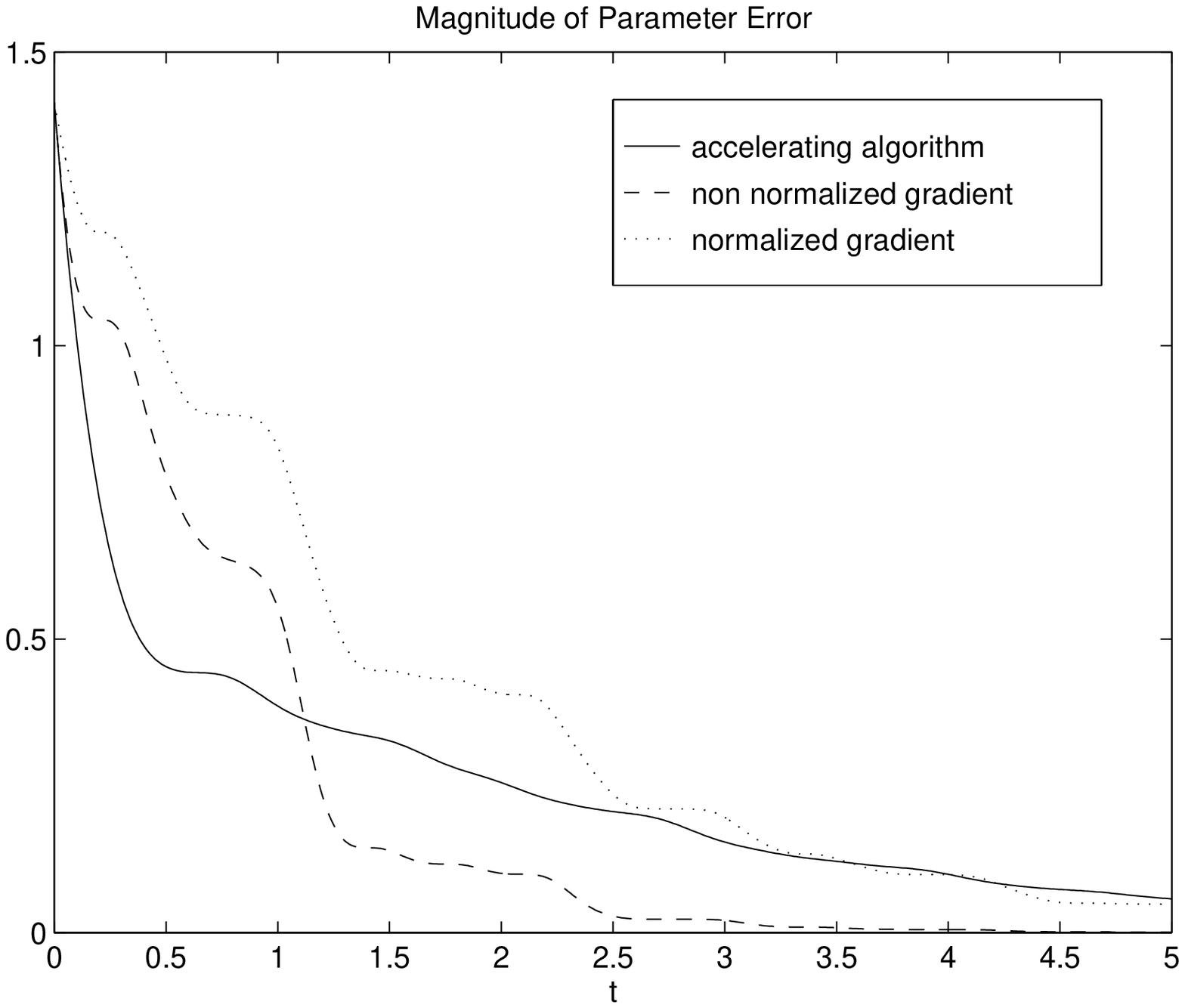} \epsfysize=2.5in \epsffile{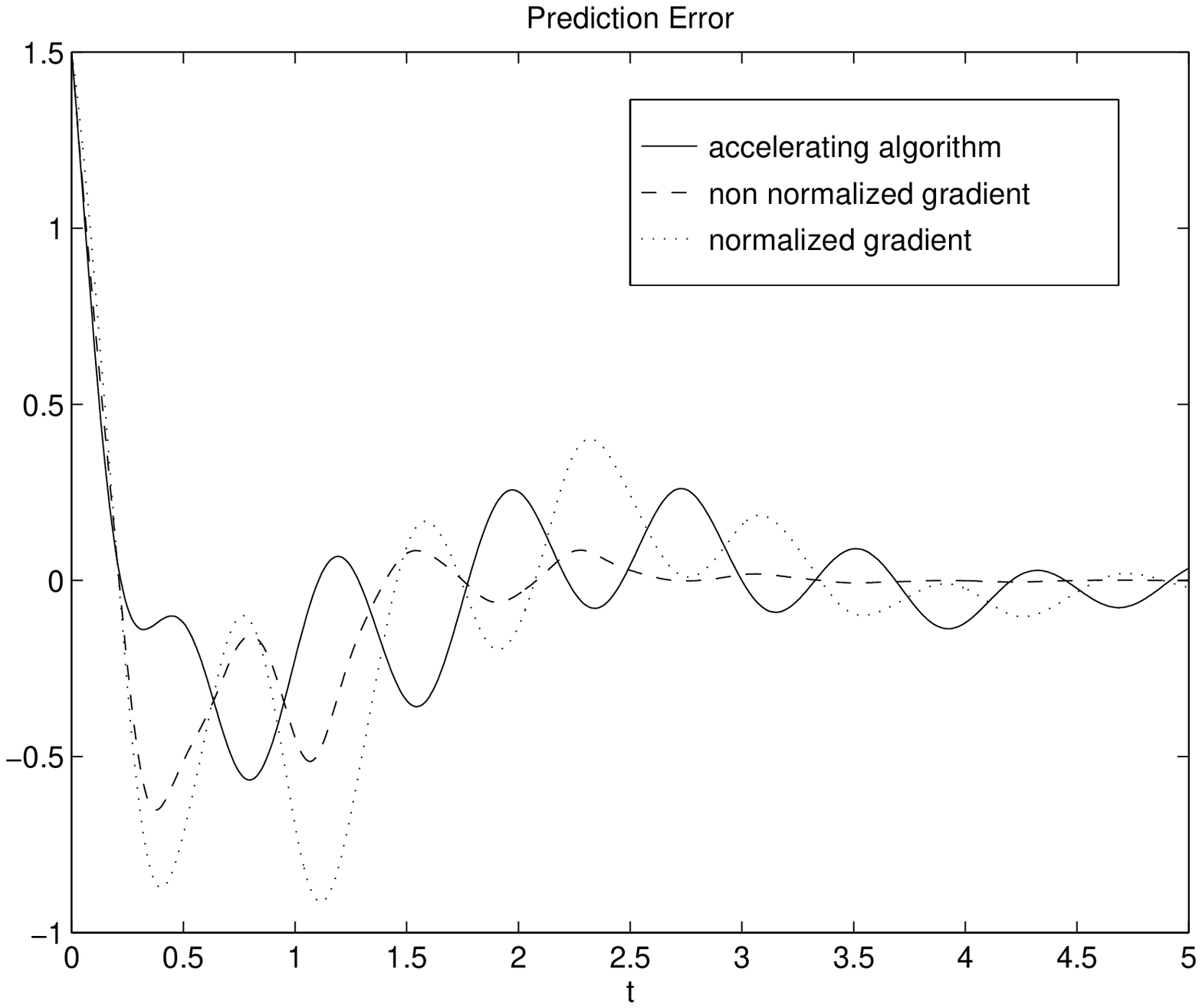}} 
 \caption{Response to sinusoidal data}
  \label{fig:sin}
\end{figure*}

\begin{figure*}[p!] 
 \epsfysize=2in
 \centerline{\epsfysize=2.5in \epsffile{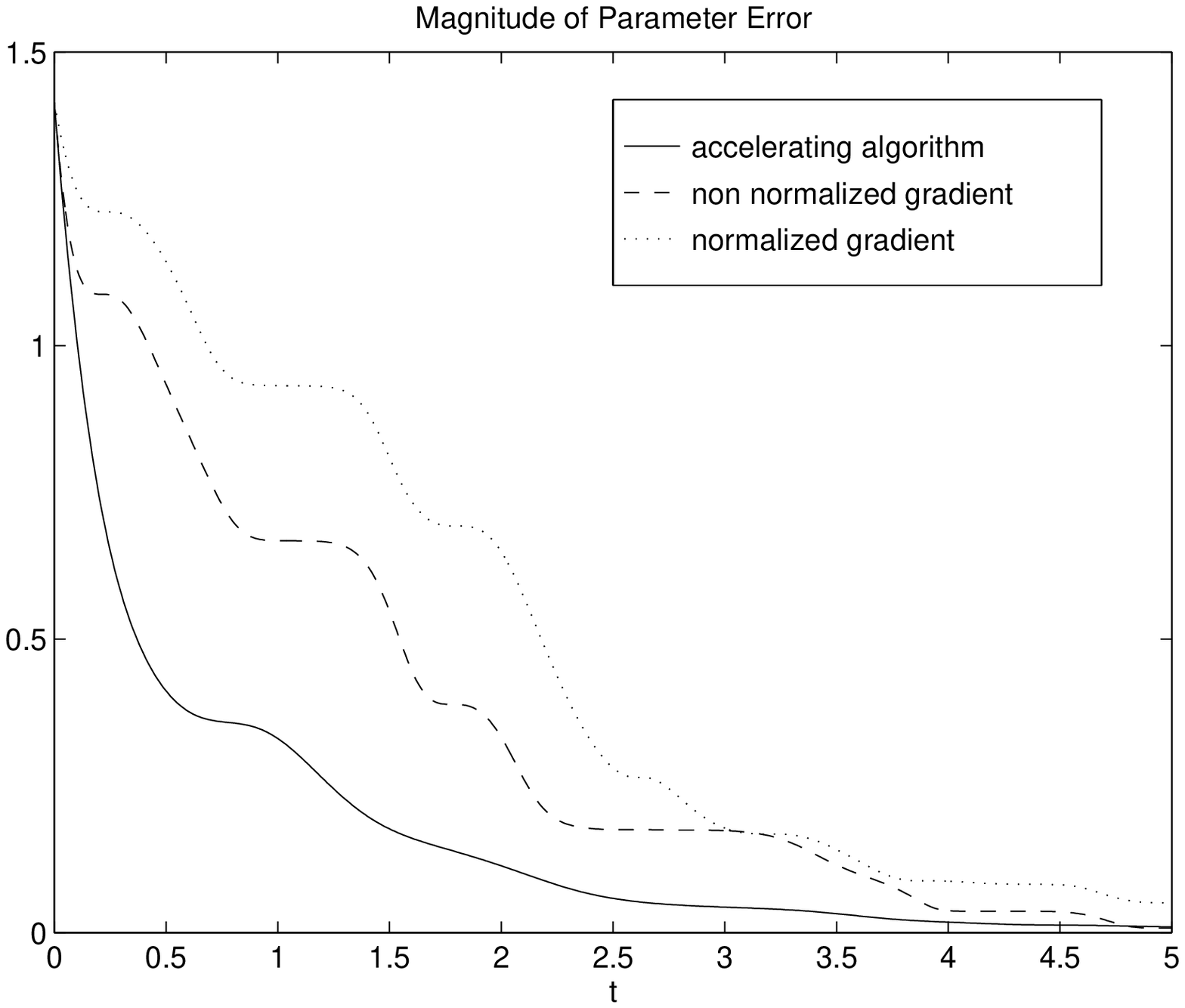} \epsfysize=2.5in \epsffile{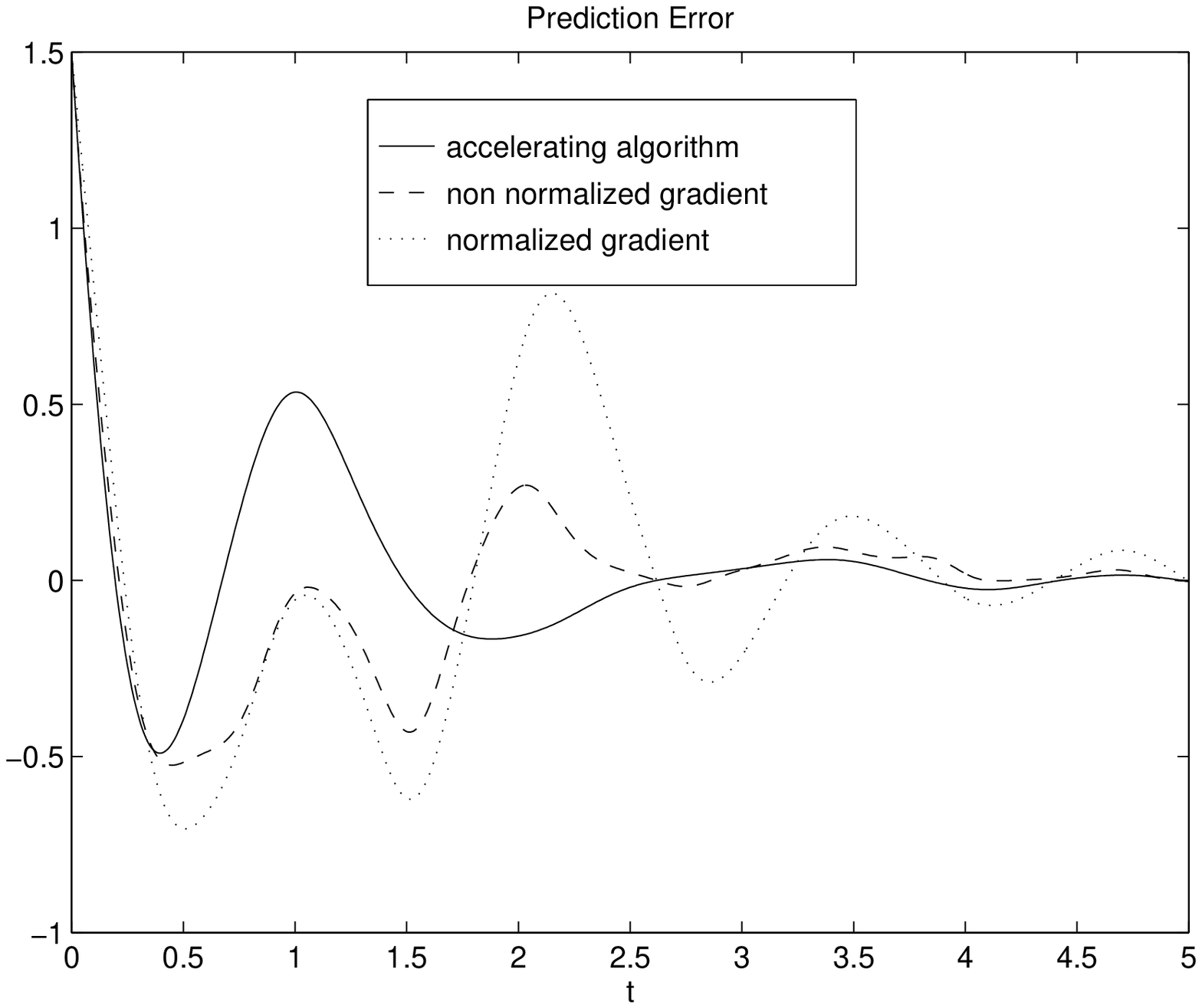}} 
 \caption{Response to exponentially increasing sinusoid}
 \label{fig:sia}
\end{figure*}

\begin{figure*}[p!] 
 \epsfysize=2in
 \centerline{\epsfysize=2.5in \epsffile{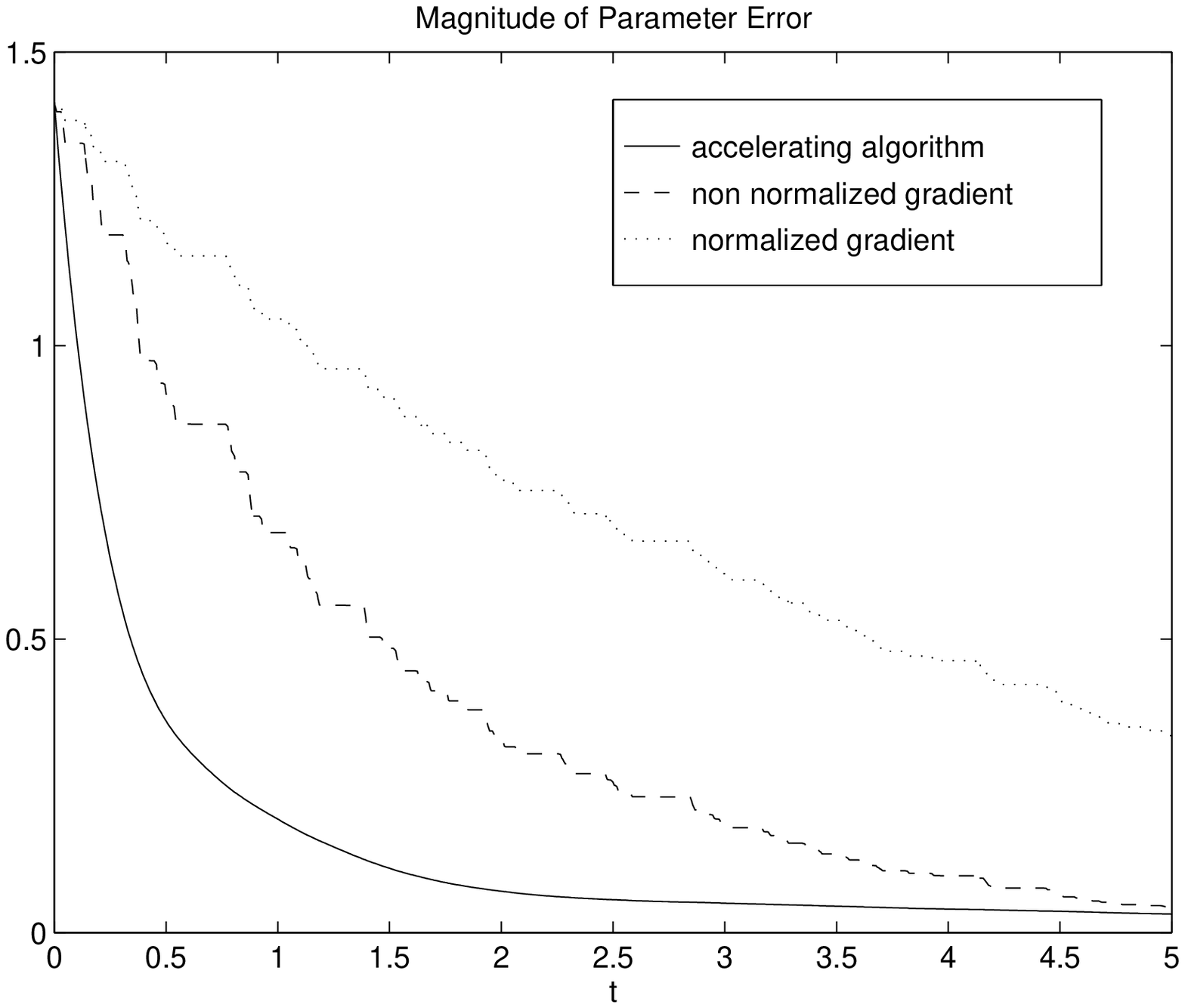} \epsfysize=2.5in \epsffile{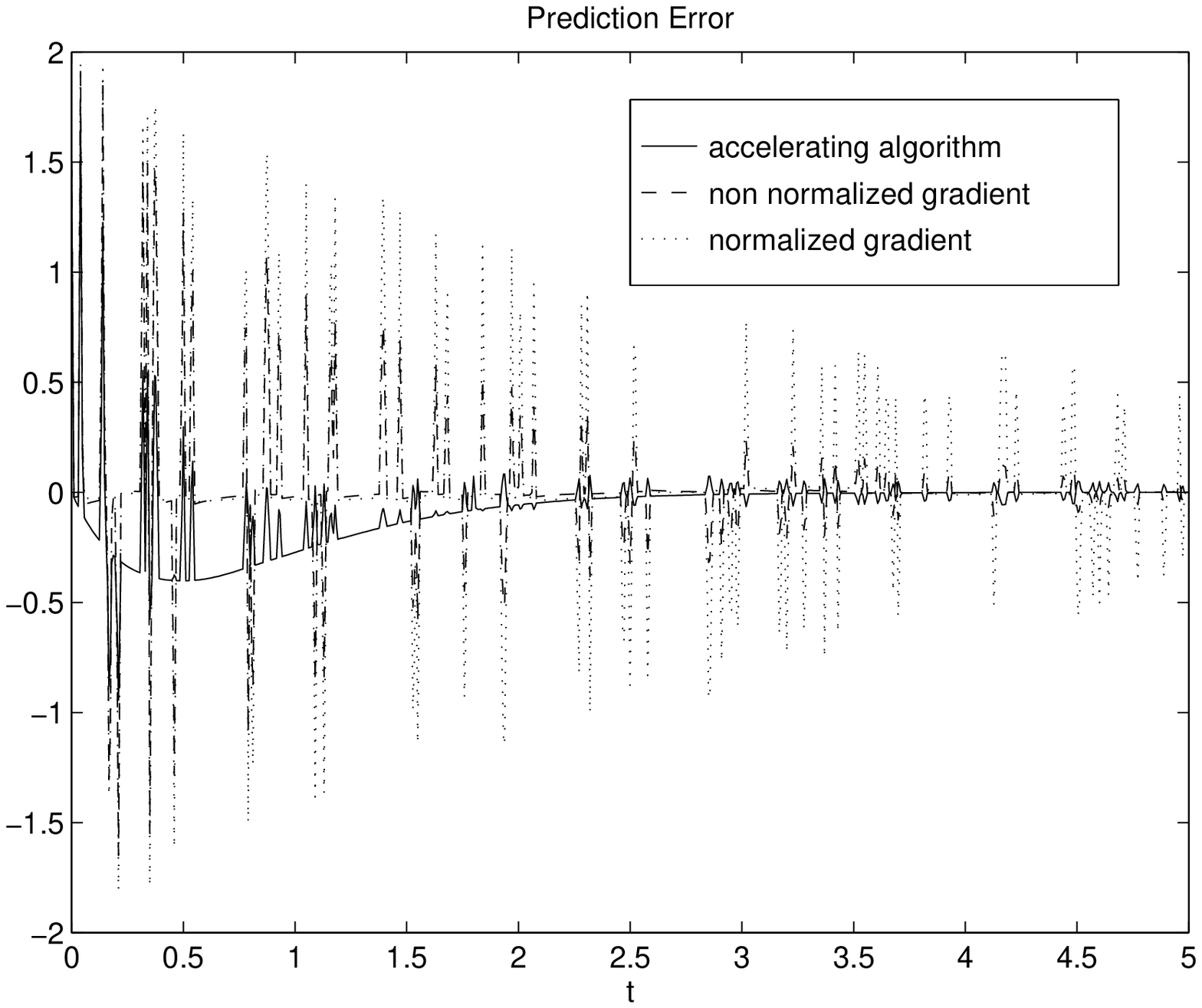}} 
 \caption{Response to pseudorandom signal}
 \label{fig:prb}
\end{figure*}

The simulations in this section compare the behavior of the accelerating tuner \rep{eq:qddot} with those of the gradient tuner \rep{eq:qdot} and of a normalized gradient one. All simulations were done in open-loop, with the regressor a two-dimensional signal, and without measurement noise. Figure~\ref{fig:step}  shows the values of $\sqrt{q^\top  q}$ and $x^\top q$ respectively when $x$ is a two-dimensional step signal. 
In Figure~\ref{fig:sin}  the regressor is a sinusoid, in 
Figure~\ref{fig:sia} an exponentially increasing sinusoid, and in 
Figure~\ref{fig:prb} a pseudorandom signal generated using \textsc{matlab}.  No effort was made to optimize the choice of gain matrices ($M_1$, $M_2$, and $M_3$ were all chosen equal to the identity), and the effect of measurement noise was not considered. The performance of the accelerating tuner is comparable, and sometimes superior, to that of the other tuners.

\section{Concluding Remarks} Other ideas related to the present one are replacing the integrator in \rep{eq:qdot} with a positive-real transfer function  \cite[page 89]{landau-book}, and using high-order tuning (\cite{morse-high-order,ortega-high-order}). High-order tuning generates as outputs $\hat{p}$ as well as its derivatives up to a given order (in this sense we might consider the present algorithm a second-order tuner), but unlike the accelerating tuner requires derivatives of $x$ up to that same order. 
We expect that accelerating tuners will find application in adaptive control of nonlinear systems and maybe in dealing with the topological incompatibility known as the ``loss of stabilizability problem'' in the adaptive control literature. 

The  stochastic analysis in \S\ref{sec:covariance} indicates that the  performance and convergence properties of the accelerating tuner, together with its moderate computational complexity, may indeed make it a desirable tool for adaptive filtering applications.
It seems that a better transient $\times$ steady-state performance compromise is achievable with the accelerating tuner than with the velocity tuner. To verify this conjecture, a study of convergence properties of the accelerating tuner and their relation with the persistence of excitation conditions is in order, as well as more extensive simulations in the presence of measurement noise. 


\paragraph{Acknowledgements}{ }The authors are grateful to \mbox{Oswaldo} L V Costa and Max Gerken for useful suggestions related to this work. The present work is an expanded version of \cite{acelera-letters}.

\bibliographystyle{plain}

\end{document}